%
\font\fifteenrm=cmr10 scaled\magstep2 
\font\fifteeni=cmmi10 scaled\magstep2
\font\fifteensy=cmsy10 scaled\magstep2
\font\fifteenbf=cmbx10 scaled\magstep2
\font\fifteentt=cmtt10 scaled\magstep2
\font\fifteenit=cmti10 scaled\magstep2
\font\fifteensl=cmsl10 scaled\magstep2
\font\fifteenam=msam10 scaled\magstep2
\font\fifteenbm=msbm10 scaled\magstep2
\font\fifteenex=cmex10 scaled\magstep2
\font\fifteensc=cmcsc10 scaled\magstep2 
\font\twelverm=cmr10 at 12pt
\font\twelvei=cmmi10 at 12pt
\font\twelvesy=cmsy10 at 12pt
\font\twelvebf=cmbx10 at 12pt
\font\twelvett=cmtt10 at 12pt
\font\twelveit=cmti10 at 12pt
\font\twelvesl=cmsl10 at 12pt
\font\twelveam=msam10 at 12pt
\font\twelvebm=msbm10 at 12pt
\font\twelveex=cmex10 at 12pt
\font\twelvesc=cmcsc10 at 12pt
\font\elevenrm=cmr10 scaled\magstephalf 
\font\eleveni=cmmi10 scaled\magstephalf
\font\elevensy=cmsy10 scaled\magstephalf
\font\elevenbf=cmbx10 scaled\magstephalf
\font\eleventt=cmtt10 scaled\magstephalf
\font\elevenit=cmti10 scaled\magstephalf
\font\elevensl=cmsl10 scaled\magstephalf
\font\elevenam=msam10 scaled\magstephalf
\font\elevenbm=msbm10 scaled\magstephalf
\font\elevenex=cmex10 scaled\magstephalf
\font\elevensc=cmcsc10 scaled\magstephalf
\font\tenrm=cmr10
\font\teni=cmmi10
\font\tensy=cmsy10
\font\tenbf=cmbx10
\font\tentt=cmtt10
\font\tenit=cmti10
\font\tensl=cmsl10
\font\tenam=msam10
\font\tenbm=msbm10
\font\tenex=cmex10
\font\tensc=cmcsc10
\font\ninerm=cmr9
\font\ninei=cmmi9
\font\ninesy=cmsy9
\font\ninebf=cmbx9
\font\ninett=cmtt9
\font\nineit=cmti9
\font\ninesl=cmsl9
\font\nineam=msam9
\font\ninebm=msbm9
\font\nineex=cmex9
\font\ninesc=cmcsc9
\font\eightrm=cmr8
\font\eighti=cmmi8
\font\eightsy=cmsy8
\font\eightbf=cmbx8
\font\eighttt=cmtt8
\font\eightit=cmti8
\font\eightsl=cmsl8
\font\eightam=msam8
\font\eightbm=msbm8
\font\eightex=cmex8
\font\eightsc=cmcsc8
\font\sevenrm=cmr7
\font\seveni=cmmi7
\font\sevensy=cmsy7
\font\sevenbf=cmbx7

\font\sevenam=msam7
\font\sevenbm=msbm7

\font\sixrm=cmr6
\font\sixi=cmmi6
\font\sixsy=cmsy6

\font\sixam=msam6
\font\sixbm=msbm6

\font\fiverm=cmr5
\font\fivei=cmmi5
\font\fivesy=cmsy5

\font\fiveam=msam5
\font\fivebm=msbm5

\font\fourrm=cmr5 at 4pt
\font\fouri=cmmi5 at 4pt
\font\foursy=cmsy5 at 4pt

\font\fouram=msam5 at 4pt
\font\fourbm=msbm5 at 4pt

\skewchar\twelvei='177 \skewchar\eleveni='177\skewchar\teni='177
\skewchar\ninei='177 \skewchar\eighti='177\skewchar\seveni='177 
\skewchar\sixi='177 \skewchar\fivei='177 \skewchar\fouri='177
\skewchar\twelvesy='60 \skewchar\elevensy='60 \skewchar\tensy='60
\skewchar\ninesy='60 \skewchar\eightsy='60 \skewchar\sevensy='60 
\skewchar\sixsy='60 \skewchar\fivesy='60 \skewchar\foursy='60
\newfam\itfam
\newfam\slfam
\newfam\bffam
\newfam\ttfam
\newfam\scfam
\newfam\amfam
\newfam\bmfam
\def\eightbig#1{{\hbox{$\left#1\vbox to 6.5pt{}\voidright $}}}
\def\eightBig#1{{\hbox{$\left#1\vbox to 7.5pt{}\voidright $}}}
\def\eightbigg#1{{\hbox{$\left#1\vbox to 10pt{}\voidright $}}}
\def\eightBigg#1{{\hbox{$\left#1\vbox to 13pt{}\voidright $}}}
\def\ninebig#1{{\hbox{$\left#1\vbox to 7.5pt{}\voidright $}}}
\def\nineBig#1{{\hbox{$\left#1\vbox to 8.5pt{}\voidright $}}}
\def\ninebigg#1{{\hbox{$\left#1\vbox to 11.5pt{}\voidright $}}}
\def\nineBigg#1{{\hbox{$\left#1\vbox to 14.5pt{}\voidright $}}}
\def\tenbig#1{{\hbox{$\left#1\vbox to 8.5pt{}\voidright $}}}
\def\tenBig#1{{\hbox{$\left#1\vbox to 9.5pt{}\voidright $}}}
\def\tenbigg#1{{\hbox{$\left#1\vbox to 12.5pt{}\voidright $}}}
\def\tenBigg#1{{\hbox{$\left#1\vbox to 16pt{}\voidright $}}}
\def\elevenbig#1{{\hbox{$\left#1\vbox to 9pt{}\voidright $}}}
\def\elevenBig#1{{\hbox{$\left#1\vbox to 10.5pt{}\voidright $}}}
\def\elevenbigg#1{{\hbox{$\left#1\vbox to 14pt{}\voidright $}}}
\def\elevenBigg#1{{\hbox{$\left#1\vbox to 17.5pt{}\voidright $}}}
\def\twelvebig#1{{\hbox{$\left#1\vbox to 10pt{}\voidright $}}}
\def\twelveBig#1{{\hbox{$\left#1\vbox to 11pt{}\voidright $}}}
\def\twelvebigg#1{{\hbox{$\left#1\vbox to 15pt{}\voidright $}}}
\def\twelveBigg#1{{\hbox{$\left#1\vbox to 19pt{}\voidright $}}}
\def\fifteenbig#1{{\hbox{$\left#1\vbox to 12pt{}\voidright $}}}
\def\fifteenBig#1{{\hbox{$\left#1\vbox to 13.5pt{}\voidright $}}}
\def\fifteenbigg#1{{\hbox{$\left#1\vbox to 18pt{}\voidright $}}}
\def\fifteenBigg#1{{\hbox{$\left#1\vbox to 23pt{}\voidright $}}}
\def\voidright{\right.\nulldelimiterspace=0pt \mathsurround=0pt }
\def\fifteenpoint{
  \textfont0=\fifteenrm \scriptfont0=\twelverm \scriptscriptfont0=\tenrm
  \def\rm{\fam0 \fifteenrm}%
  \textfont1=\fifteeni \scriptfont1=\twelvei \scriptscriptfont1=\teni
  \textfont2=\fifteensy \scriptfont2=\twelvesy \scriptscriptfont2=\tensy
  \textfont3=\fifteenex \scriptfont3=\fifteenex \scriptscriptfont3=\fifteenex
  \def\it{\fam\itfam\fifteenit}\textfont\itfam=\fifteenit
  \def\sl{\fam\slfam\fifteensl}\textfont\slfam=\fifteensl
  \def\bf{\fam\bffam\fifteenbf}\textfont\bffam=\fifteenbf 
    \scriptfont\bffam=\twelvebf\scriptscriptfont\bffam=\tenbf
  \def\tt{\fam\ttfam\fifteentt}\textfont\ttfam=\fifteentt
  \def\sc{\fam\scfam\fifteensc}\textfont\scfam=\fifteensc
  \def\am{\fam\amfam\fifteenam}\textfont\amfam=\fifteenam
    \scriptfont\amfam=\twelveam\scriptscriptfont\amfam=\tenam
  \def\bm{\fam\bmfam\fifteenbm}\textfont\bmfam=\fifteenbm
    \scriptfont\bmfam=\twelvebm\scriptscriptfont\bmfam=\tenbm
  \baselineskip=21pt \rm
  \let\big=\fifteenbig\let\Big=\fifteenBig\let\bigg=\fifteenbigg
  \let\Bigg=\fifteenBigg}
\def\twelvepoint{
  \textfont0=\twelverm \scriptfont0=\ninerm \scriptscriptfont0=\sevenrm
  \def\rm{\fam0 \twelverm}%
  \textfont1=\twelvei \scriptfont1=\ninei \scriptscriptfont1=\seveni
  \textfont2=\twelvesy \scriptfont2=\ninesy \scriptscriptfont2=\sevensy
  \textfont3=\twelveex \scriptfont3=\twelveex \scriptscriptfont3=\twelveex
  \def\it{\fam\itfam\twelveit}\textfont\itfam=\twelveit
  \def\sl{\fam\slfam\twelvesl}\textfont\slfam=\twelvesl
  \def\bf{\fam\bffam\twelvebf}\textfont\bffam=\twelvebf 
    \scriptfont\bffam=\ninebf\scriptscriptfont\bffam=\sevenbf
  \def\tt{\fam\ttfam\twelvett}\textfont\ttfam=\twelvett
  \def\sc{\fam\scfam\twelvesc}\textfont\scfam=\twelvesc
  \def\am{\fam\amfam\twelveam}\textfont\amfam=\twelveam
    \scriptfont\amfam=\nineam\scriptscriptfont\amfam=\sevenam
  \def\bm{\fam\bmfam\twelvebm}\textfont\bmfam=\twelvebm
    \scriptfont\bmfam=\ninebm\scriptscriptfont\bmfam=\sevenbm
  \baselineskip=17.8pt \rm 
  \def\looselineskip{\baselineskip=18.5pt plus 1.8pt}%
  \def\tightlineskip{\baselineskip=16.5pt}%
  \def\verytightlineskip{\baselineskip=15pt}%
  \let\big=\twelvebig\let\Big=\twelveBig\let\bigg=\twelvebigg
  \let\Bigg=\twelveBigg  }
\def\elevenpoint{
  \textfont0=\elevenrm \scriptfont0=\ninerm \scriptscriptfont0=\sixrm
  \def\rm{\fam0 \elevenrm}%
  \textfont1=\eleveni \scriptfont1=\ninei \scriptscriptfont1=\sixi
  \textfont2=\elevensy \scriptfont2=\ninesy \scriptfont2=\sixsy 
  \textfont3=\elevenex \scriptfont3=\elevenex \scriptfont3=\elevenex
  \def\it{\fam\itfam\elevenit}\textfont\itfam=\elevenit
  \def\sl{\fam\slfam\elevensl}\textfont\slfam=\elevensl
  \def\bf{\fam\bffam\elevenbf}\textfont\bffam=\elevenbf
  \def\tt{\fam\ttfam\eleventt}\textfont\ttfam=\eleventt
  \def\sc{\fam\scfam\elevensc}\textfont\scfam=\elevensc
  \def\am{\fam\amfam\elevenam}\textfont\amfam=\elevenam
    \scriptfont\amfam=\nineam\scriptscriptfont\amfam=\sixam
  \def\bm{\fam\bmfam\elevenbm}\textfont\bmfam=\elevenbm
    \scriptfont\bmfam=\ninebm\scriptscriptfont\bmfam=\sixbm
  \baselineskip=15.1pt \rm
  \def\looselineskip{\baselineskip=16pt plus 1.5pt}%
  \def\tightlineskip{\baselineskip=14pt}%
  \def\verytightlineskip{\baselineskip=13pt}%
  \let\big=\elevenbig\let\Big=\elevenBig\let\bigg=\elevenbigg
  \let\Bigg=\elevenBigg  }
\def\tenpoint{
  \textfont0=\tenrm \scriptfont0=\eightrm \scriptscriptfont0=\fiverm
  \def\rm{\fam0 \tenrm}%
  \textfont1=\teni \scriptfont1=\eighti \scriptscriptfont1=\fivei
  \textfont2=\tensy \scriptfont2=\eightsy \scriptfont2=\fivesy 
  \textfont3=\tenex \scriptfont3=\tenex \scriptfont3=\tenex
  \def\it{\fam\itfam\tenit}\textfont\itfam=\tenit
  \def\sl{\fam\slfam\tensl}\textfont\slfam=\tensl
  \def\bf{\fam\bffam\tenbf}\textfont\bffam=\tenbf
  \def\tt{\fam\ttfam\tentt}\textfont\ttfam=\tentt
  \def\sc{\fam\scfam\tensc}\textfont\scfam=\tensc
  \def\am{\fam\amfam\tenam}\textfont\amfam=\tenam
    \scriptfont\amfam=\eightam \scriptscriptfont\amfam=\fiveam
  \def\bm{\fam\bmfam\tenbm}\textfont\bmfam=\tenbm
    \scriptfont\bmfam=\eightbm \scriptscriptfont\bmfam=\fivebm
  \baselineskip=14pt \rm
  \def\looselineskip{\baselineskip=14.8pt plus1.5pt}
  \def\tightlineskip{\baselineskip=13.6pt}%
  \def\verytightlineskip{\baselineskip=13pt}%
  \let\big=\tenbig\let\Big=\tenBig\let\bigg=\tenbigg\let\Bigg=\tenBigg  }
\def\ninepoint{
  \textfont0=\ninerm \scriptfont0=\sevenrm \scriptscriptfont0=\fourrm
  \def\rm{\fam0 \ninerm}%
  \textfont1=\ninei \scriptfont1=\seveni \scriptscriptfont1=\fouri
  \textfont2=\ninesy \scriptfont2=\sevensy \scriptfont2=\foursy 
  \textfont3=\nineex \scriptfont3=\nineex \scriptfont3=\nineex
  \def\it{\fam\itfam\nineit}\textfont\itfam=\nineit
  \def\sl{\fam\slfam\ninesl}\textfont\slfam=\ninesl
  \def\bf{\fam\bffam\ninebf}\textfont\bffam=\ninebf
  \def\tt{\fam\ttfam\ninett}\textfont\ttfam=\ninett
  \def\sc{\fam\scfam\ninesc}\textfont\scfam=\ninesc
  \def\am{\fam\amfam\nineam}\textfont\amfam=\nineam
    \scriptfont\amfam=\nineam\scriptscriptfont\amfam=\fouram
  \def\bm{\fam\bmfam\ninebm}\textfont\bmfam=\ninebm
    \scriptfont\bmfam=\ninebm\scriptscriptfont\bmfam=\fourbm
  \baselineskip=12.6pt \rm
  \let\big=\ninebig\let\Big=\nineBig\let\bigg=\ninebigg
  \let\Bigg=\nineBigg  }
\def\eightpoint{
  \textfont0=\eightrm \scriptfont0=\fiverm \scriptscriptfont0=\fourrm
  \def\rm{\fam0 \eightrm}%
  \textfont1=\eighti \scriptfont1=\fivei \scriptscriptfont1=\fouri
  \textfont2=\eightsy \scriptfont2=\fivesy \scriptfont2=\foursy 
  \textfont3=\eightex \scriptfont3=\eightex \scriptfont3=\eightex
  \def\it{\fam\itfam\eightit}\textfont\itfam=\eightit
  \def\sl{\fam\slfam\eightsl}\textfont\slfam=\eightsl
  \def\bf{\fam\bffam\eightbf}\textfont\bffam=\eightbf
  \def\tt{\fam\ttfam\eighttt}\textfont\ttfam=\eighttt
  \def\sc{\fam\scfam\eightsc}\textfont\scfam=\eightsc
  \def\am{\fam\amfam\eightam}\textfont\amfam=\eightam
    \scriptfont\amfam=\eightam\scriptscriptfont\amfam=\fouram
  \def\bm{\fam\bmfam\eightbm}\textfont\bmfam=\eightbm
    \scriptfont\bmfam=\eightbm\scriptscriptfont\bmfam=\fourbm
  \baselineskip=11.2pt \rm
  \let\big=\eightbig\let\Big=\eightBig\let\bigg=\eightbigg
  \let\Bigg=\eightBigg  }

\twelvepoint
\nopagenumbers
\hsize=6in\vsize=8.8in

\parskip=1pt plus 1pt

\newif\ifSpecialhead\Specialheadfalse
\newbox\specialheadbox

\def\specialhead #1\par{\Specialheadtrue\setbox\specialheadbox=\hbox{#1}}
\headline={{\ifSpecialhead\box\specialheadbox\global\Specialheadfalse\else
     \ifnum\pageno<0{\hfill\quad{\twelvebf\folio}}%
     \else\ifnum\pageno<2\hfill
     \else\hfill\twelvepoint\sc\firstmark\quad{\twelvebf\folio}\fi\fi\fi}}

\def\title#1\par{\bigskip{\def\cr{\par\center}\center\fifteenbf #1\par}\medskip}
\def\subtitle#1\par{\centerline{\fifteenrm #1}\medskip}
\def\author#1\par{\medskip{\def\cr{\par\center\twelvesc}\fifteensc\center#1\par}}
\def\center#1\par{\hfil #1\hfil\par}
\def\abstract.#1\par{\message{Abstract.}%
                    \medskip{\narrower\narrower\tenpoint\tightlineskip
                        \noindent{\bf Abstract.}#1\par}\medskip\noindent}
\def\bigabstract.#1\par{\message{Abstract.}%
                         \medskip{\narrower\narrower\tightlineskip
                         \noindent{\bf Abstract. }#1\par}\medskip\noindent}
\def\acknowledgement#1\par{\footnote{}{#1}}
\def\sectionskip{\Goodbreak\vskip 25pt plus 15pt minus 5pt}
\def\secnumber{\ifquiet
               \else\ifNoSections
                    \else\sectionsymbol\the\secno\quad\fi\fi}
\def\section#1\par{ \NoSectionsfalse\par\sectionskip\proofdepth=0\claimno=0
 \ifquiet\else\advance\secno by1\fi\toks0={#1}
 \immediate\write16{\ifquiet\else Section \the\secno\space\fi
                    \the\toks0}%
 \mark{\secnumber #1}%
 {\fifteenpoint\bf\noindent\secnumber #1}\nobreak\bigskip\quietoff
 \nobreak\noindent}
\def\quiet{\quiettrue}

\def\quietoff{\ifQUIET\else\quietfalse\fi}
\newif\ifquiet
\newif\ifQUIET
\newif\ifNoSections
\newcount\claimtype
\newcount\secno
\newcount\claimno
\newcount\subclaimno
\newcount\subsubclaimno
\newcount\subsubsubclaimno
\newcount\proofdepth
\def\subclaimnumber{\ifquiet\else\ifcase\subclaimno\or A\or B\or C\or D\or E\or
     F\or G\or H\or I\or J\or K\or L\or M\or N\or O\or P\fi\fi}
\def\subsubclaimnumber{\ifquiet\else\ifcase\subsubclaimno\or i\or ii\or iii\or 
   iv\or v\or vi\or vii\or viii\or ix\or x\or xi\or xii\or xiii\or xiv\fi\fi}
\def\subsubsubclaimnumber{\ifquiet\else\ifcase\subsubsubclaimno\or a\or b\or 
   c\or d\or e\or f\or g\or \or h\or i\or j\or k\or l\or m\or n\or o\fi\fi}
\def\claimtag{\ifquiet\else
  \ifNoSections
    \ifcase\proofdepth\the\claimno%
    \or\the\claimno.\subclaimnumber
    \or\the\claimno.\subclaimnumber.\subsubclaimnumber
    \or\the\claimno.\subclaimnumber.\subsubclaimnumber
                                                .\subsubsubclaimnumber\fi
  \else
    \ifcase\proofdepth\the\secno.\the\claimno
    \or\the\secno.\the\claimno.\subclaimnumber
    \or\the\secno.\the\claimno.\subclaimnumber.\subsubclaimnumber
    \or\the\secno.\the\claimno.\subclaimnumber.\subsubclaimnumber
                                                .\subsubsubclaimnumber\fi\fi\fi}
\secno=0\claimno=0\proofdepth=0\subclaimno=0\subsubclaimno=0\subsubsubclaimno=0
\NoSectionstrue
\newbox\qedbox
\def\claimname{\ifcase\claimtype Theorem\or Lemma\or Claim\or Corollary\or
               Question\or Definition\or Remark\or Conjecture\fi}
\def\preclaimskip{\removelastskip
    \ifcase\claimtype\goodbreak\vskip 8pt plus 4pt minus 2pt
                  \or\goodbreak\vskip 6pt plus 4pt minus 1pt
                  \or\goodbreak\vskip 5pt plus 4pt minus 1pt
                  \or\goodbreak\vskip 8pt plus 4pt minus 2pt
                  \or\vskip 7pt plus 4pt minus 2pt
                  \or\vskip 7pt plus 4pt minus 2pt
                  \or\vskip 7pt plus 4pt minus 2pt
                  \or\goodbreak\vskip 8pt plus 4pt minus 2pt\fi}
\def\postclaimskip{\ifcase\claimtype         \vskip 4pt plus 2pt minus 2pt
                                          \or\vskip 3pt plus 2pt minus 2pt
                                          \or\vskip 2pt plus 2pt minus 1pt
                                          \or\vskip 4pt plus 2pt minus 2pt
                                          \or\vskip 1pt plus 2pt 
                                          \or\vskip 4pt plus 4pt 
                                          \or\vskip 3pt plus 2pt
                                          \or\vskip 4pt plus 2pt minus 2pt\fi}
\def\claimfont{\ifcase\claimtype
                  \sl\or\sl\or\sl\or\sl\or\sl\or\rm\or\rm\or\sl\fi}
\def\advancetag{\ifcase\proofdepth\advance\claimno by1
                               \or\advance\subclaimno by1
                               \or\advance\subsubclaimno by1
                               \or\advance\subsubsubclaimno by1\fi}
\def\sayclaim#1.#2 #3\par{\ifquiet\else\advancetag\fi
    \preclaimskip\setbox1=\hbox{#1}\setbox2=\hbox{#2}%
    \toks0={#1 }
    \immediate\write16{\ifdim\wd1>0pt\the\toks0
                       \else\claimname\space\fi \claimtag.}%
    \vbox{\noindent
    {\bf\ifdim\wd1=0pt \claimname\else #1\fi\ifquiet.\else\ \claimtag{\ifNoSections.\fi}\fi}%
    \enspace{\ifdim\wd2>0pt\sc #2\enspace\fi}%
    {\claimfont #3\par}}\postclaimskip\quietoff}
\def\theorem{\claimtype=0\sayclaim}
\def\lemma{\claimtype=1\sayclaim}

\def\corollary{\claimtype=3\sayclaim}
\def\question{\claimtype=4\sayclaim}

\def\remark{\claimtype=6\sayclaim}

\def\point#1. #2\par{\item{\rm #1.}#2\par}
\def\points#1\cr\par{\medskip\vbox{\let\cr=\point\point#1\par}\par}
\def\df{\it}
\def\prooffont{}
\def\proofsize{}
\def\proofindent{}
\def\proofskip{\badbreak\ifcase\claimtype    \vskip 3pt plus 2pt minus 2pt
                                          \or\vskip 2pt plus 2pt minus 2pt
                                          \or\vskip 1pt plus 2pt minus 1pt
                                          \or\vskip 3pt plus 2pt minus 2pt
                                          \or\vskip 1pt plus 2pt 
                                          \or\vskip 2pt plus 4pt 
                                          \or\vskip 1pt plus 2pt
                                          \or\vskip 3pt plus 2pt minus 2pt\fi}

\def\Goodbreak{\vskip0pt plus.5in\penalty-1000\vskip0pt plus-.5in}
\def\goodbreak{\penalty-500}
\def\badbreak{\penalty500}
\def\Badbreak{\penalty1000}
\def\proof{\message{proof}\removelastskip\Badbreak\proofskip\begingroup
  \advance\proofdepth by1
  \setbox\qedbox=\hbox{\halmos\raise2pt\hbox{\fiverm\claimname}}%
  \prooffont\proofsize\proofindent\noindent{\bf Proof: }}
\def\proofof#1:{\message{proof}\removelastskip\Badbreak\proofskip\begingroup
  \advance\proofdepth by1
  \setbox\qedbox=\hbox{\halmos\raise2pt\hbox{\fiverm#1}}%
  \prooffont\proofsize\proofindent\noindent{\bf Proof of #1: }}
\def\cite[#1]{[{\tenrm{#1}}]\message{[#1]}}
\edef\ref#1{\expandafter\global\expandafter\edef#1{\noexpand\claimtag}}
\newwrite\notes
\openout\notes=\jobname.notes
\long\def\unexpandedwrite#1#2{\def\finwrite{\write#1}%
   {\aftergroup\finwrite\aftergroup{\sanitize#2\endsanity}}}
\def\sanitize{\futurelet\next\sanswitch}
\let\stoken=\space
\def\sanswitch{\ifx\next\endsanity
  \else\ifcat\noexpand\next\stoken\aftergroup\space\let\next=\eat
   \else\ifcat\noexpand\next\bgroup\aftergroup{\let\next=\eat
    \else\ifcat\noexpand\next\egroup\aftergroup}\let\next=\eat
     \else\let\next=\copytoken\fi\fi\fi\fi \next}
\def\eat{\afterassignment\sanitize \let\next= }
\long\def\copytoken#1{\ifcat\noexpand#1\relax\aftergroup\noexpand
  \else\ifcat\noexpand#1\noexpand~\aftergroup\noexpand\fi\fi
  \aftergroup#1\sanitize}
\def\endsanity\endsanity{}

\def\note#1#2{\hbox to2in{\strut#1\quad\dotfill\quad#2}}
\def\boxit#1{\setbox4=\hbox{\kern1pt#1\kern1pt}
  \hbox{\vrule\vbox{\hrule\kern1pt\box4\kern1pt\hrule}\vrule}}
\def\halmos{\hbox{\am\char'3}} 
\def\qed#1\par{\message{.                                }\setbox1=\hbox{#1}%
  \ifdim\wd1>0pt\setbox\qedbox=\hbox{\halmos\raise2pt\hbox{\fiverm #1}}\fi
  \kern5pt\lower 2pt\hbox{\box\qedbox}\proofskip\goodbreak\endgroup}

\def\sectionsymbol{\S}

\def\g{\gamma}
\def\a{\alpha}
\def\b{\beta}
\def\d{\delta}

\def\l{\lambda}
\def\z{\zeta}
\def\I1{\mathop{\hbox{\sc i}_1}}
\def\w{\omega}

\def\R{{\mathchoice{\hbox{\bm R}}{\hbox{\bm R}}
         {\hbox{\tenbm R}}{\hbox{\sevenbm R}}}}

\def\elesub{\prec}

\def\unifto{\buildrel\lower 7pt\hbox{$\to$}\over\to}

\def\from{\mathbin{\vbox{\baselineskip=3pt\lineskiplimit=0pt
                         \hbox{.}\hbox{.}\hbox{.}}}}

\def\in{\mathrel{\mathchoice{\raise 
1pt\hbox{$\scriptstyle\cal\char'62$}}
         {\raise 1pt\hbox{$\scriptstyle\cal\char'62$}}
         {\raise .5pt\hbox{$\scriptscriptstyle\cal\char'62$}}
         {\hbox{$\scriptscriptstyle\cal\char'62$}}}\penalty700{}}
\def\ni{\mathrel{\mathchoice{\raise 1pt\hbox{$\scriptstyle\cal\char'63$}}
                   {\raise 1pt\hbox{$\scriptstyle\cal\char'63$}}
                   {\raise .5pt\hbox{$\scriptscriptstyle\cal\char'63$}}
                   {\hbox{$\scriptscriptstyle\cal\char'63$}}}\penalty700}
\def\of{\mathrel{\mathchoice{\raise 1pt\hbox{$\scriptstyle\subseteq$}}
                   {\raise 1pt\hbox{$\scriptstyle\subseteq$}}
                   {\raise .5pt\hbox{$\scriptscriptstyle\subseteq$}}
                   {\hbox{$\scriptscriptstyle\subseteq$}}}}
\def\fo{\mathrel{\mathchoice{\raise 1pt\hbox{$\scriptstyle\supseteq$}}
                   {\raise 1pt\hbox{$\scriptstyle\supseteq$}}
                   {\raise .5pt\hbox{$\scriptscriptstyle\supseteq$}}
                   {\hbox{$\scriptscriptstyle\supseteq$}}}}
\def\notin{\mathrel{\mathchoice
  {\raise 1pt\hbox{\rlap{$\scriptstyle\;|$}$\scriptstyle\cal\char'62$}}
  {\raise 1pt\hbox{\rlap{$\scriptstyle\kern2pt 
          |$}$\scriptstyle\cal\char'62$}}
  {\raise .5pt\hbox{\rlap{$\scriptscriptstyle\, |$}$\scriptscriptstyle
      \cal\char'62$}}
  {\hbox{\rlap{$\scriptscriptstyle\, |$}$\scriptscriptstyle
     \cal\char'62$}}}%
  \penalty700}

\def\and{\mathrel{\kern1pt\&\kern1pt}}

\def\cross{\times}

\def\lte{\mathrel{\scriptstyle\leq}}
\def\tlt{\triangleleft}

\def\[#1]{\left[\vphantom{\bigm|}#1\right]}
\def\<#1>{\langle\,#1\,\rangle}

\def\sat{\models}

\def\restrict{\mathbin{\mathchoice{\hbox{\am\char'26}}{\hbox{\am\char'26}}{\hbox{\eightam\char'26}}{\hbox{\sixam\char'26}}}}

\def\st{\mid}
\def\seq<#1>{{\def\st{\mid\penalty650}\left<\,#1\,\right>}}

\def\set#1{\{\,#1\,\}}

\def\th{{\hbox{\fiverm th}}}

\def\lttheta{{\raise 1pt\hbox{$\scriptstyle<$}\theta}}

\def\I1{\mathop{\hbox{\sc i}_1}}

\def\jump{\triangle}
\def\Jump{\fulltriangle}

\def\jump{{\!\hbox{\nineam\char'117}}}
\def\Jump{{\!\hbox{\nineam\char'110}}}

\font\arrow=line10 scaled \magstep1
\def\makeline#1.{\hbox{\arrow\char#1}}
\def\makearrow#1.#2.{\hbox{\arrow\char#1\llap{\char#2}}}
\def\definelinesandarrows#1.#2.#3.#4.#5.{
   \expandafter\edef\csname#4line\endcsname{\makeline#1.}
   \expandafter\edef\csname#4arrow\endcsname{\makearrow#1.#2.}
   \expandafter\edef\csname#5line\endcsname{\makeline#1.}
   \expandafter\edef\csname#5arrow\endcsname{\makearrow#1.#3.}}
\definelinesandarrows 0.18.9.ne.sw.
\definelinesandarrows 1.21.11.nnne.sssw.
\definelinesandarrows 2.14.13.nnnne.ssssw.
\definelinesandarrows 3.23.15.nnnnne.sssssw.
\definelinesandarrows 4.23.15.nnnnnne.ssssssw.
\definelinesandarrows 10.30.29.nne.ssw.
\definelinesandarrows 16.49.41.neeeeee.swwwwww.
\definelinesandarrows 17.51.43.neeee.swwww.
\definelinesandarrows 19.55.47.nehuh.swhuh.
\definelinesandarrows 24.58.41.neeeeeee.swwwwwww.
\definelinesandarrows 26.62.9.neee.swww.
\definelinesandarrows 33.49.25.neeeee.swwwww.
\definelinesandarrows 35.62.61.nee.sww.
\definelinesandarrows 64.82.73.se.nw.
\definelinesandarrows 65.85.75.ssse.nnnw.
\definelinesandarrows 66.78.77.sssse.nnnnw.
\definelinesandarrows 67.87.79.ssssse.nnnnnw.
\definelinesandarrows 68.87.79.sssssse.nnnnnnw.
\definelinesandarrows 74.94.93.sse.nnw.
\definelinesandarrows 80.113.105.seeeeee.nwwwwww.
\definelinesandarrows 81.115.107.seeee.nwwww.
\definelinesandarrows 99.126.125.see.nww.
\def\sejoin#1#2{\setbox1=\hbox{#1}\setbox2=\hbox{#2}%
  \hbox{\vbox{\hbox{\copy1\kern\wd2}\nointerlineskip
              \hbox{\kern\wd1\box2}}}}
\def\nejoin#1#2{\setbox1=\hbox{#1}\setbox2=\hbox{#2}%
  \hbox{\vbox{\hbox{\kern\wd1\copy2}\nointerlineskip\hbox{\copy1\kern\wd2}}}}
\newdimen\hnudge
\newdimen\vnudge
\newdimen\hnudgedefault
\newdimen\vnudgedefault

\def\SEdefaultnudge{\hnudge=-16pt\vnudge=20pt}
\def\Edefaultnudge{\hnudge=-25pt\vnudge=6pt}

\def\longEdefaultnudge{\hnudge=-5pt\vnudge=6pt}
\def\nudgeright#1pt{\advance\hnudge by#1pt}
\def\nudgeleft#1pt{\advance\hnudge by-#1pt}
\def\nudgeup#1pt{\advance\vnudge by#1pt}
\def\nudgedown#1pt{\advance\vnudge by-#1pt}
\def\label#1{\smash{\llap{\kern\hnudge
                   \raise\vnudge\rlap{$\scriptstyle#1$}\hfill}}}

\def\SEarrow{\SEdefaultnudge
             \sejoin\seeline{\sejoin\seeline{\sejoin\seeline\seearrow}}}

\def\Earrow{\Edefaultnudge\setbox1=\hbox{\SEarrow}
 \hbox{\raise 2pt\hbox{\vrule height-.4pt depth.8ptwidth\wd1\kern2pt
       \llap{\arrow\char'55}}}}
\def\longEarrow{\longEdefaultnudge\setbox1=\hbox{\SEarrow}
      \rlap{\hskip-1.25\wd1\raise 2pt
            \hbox{\vrule height-.4pt depth.8ptwidth2.5\wd1\kern2pt
            \llap{\arrow\char'55}}}}

\def\ilt{<_{\scriptscriptstyle\infty}}
\def\ilte{\lte_{\scriptscriptstyle\infty}}
\def\iequiv{\equiv_{\scriptscriptstyle\infty}}
\def\halts{\downarrow}

\center [submitted to the Archive for Mathematical Logic]

\title Post's problem for supertasks has both\cr
       positive and negative solutions\cr

\author Joel David Hamkins\cr
      Kobe University and the\cr
      City University of New York\cr
      {\ninett http://www.library.csi.cuny.edu/users/hamkins}\cr

\author Andrew Lewis\cr
      Virginia Commonwealth University\cr
      {\ninett http://saturn.vcu.edu/$\tilde{\phantom{a}}$amlewis/}\cr

\abstract. The infinite time Turing machine analogue of Post's problem, the question whether there are supertask degrees between $0$ and the supertask jump $0^\jump$, has in a sense both positive and negative solutions. Namely, in the context of the reals there are no degrees between $0$ and $0^\jump$, but in the context of {\it sets} of reals, there are; indeed, there are incomparable semi-decidable supertask degrees. Both arguments employ a kind of transfinite-injury construction which generalizes canonically to oracles.

\acknowledgement The first author wishes to thank his hosts at Kobe University and the Japan Society for the Promotion of Science for their gracious support during the 1998 calendar year, as well as the CUNY Scholar Incentive program. His research has been supported in part by a grant from the PSC-CUNY Research Foundation.

The supertask model of infinitary computation extends the operation of an ordinary Turing machine into transfinite ordinal time. The theory, of which we gave a full introductory account in \cite[HamLew], leads naturally to notions of computability and semi-decidability for reals and sets of reals, as well as an oracle concept and two jump operators: the lightface $\jump$ and boldface $\Jump$ jumps. Left open in \cite[HamLew] was the supertask analogue of Post's problem, the question whether there are semi-decidable degrees between $0$ and the jump $0^\jump$, or indeed, whether there are any degrees between $0$ and $0^\jump$ at all. Here, we prove that the answer is in a sense both positive and negative. 

\quiet\theorem Main Theorem. Depending on the context, the supertask analogue of Post's problem has both positive and negative solutions. Specifically:
\points 1. In the context of the reals, there are no degrees strictly between $0$ and $0^\jump$.\cr
            2. In the context of sets of reals, there are degrees strictly between $0$ and $0^\jump$. Indeed, there are incomparable semi-decidable degrees.\cr

While the first statement illustrates how sharply the supertask recursion theory departs from the classical theory, the second hints at a broader analogy in the context of degrees on sets of reals. Indeed, the second statement is proved by adapting the classical Friedburg-Munchnik priority argument to construct countable incomparable semi-decidable sets. Both of these arguments employ a kind of transfinite-injury construction which relativizes canonically to any oracle $A$, so that first, there is no real $z$ such that $A\oplus z$ is strictly between $A$ and $A^\jump$, and second, there are $A$-incomparable $A$-semidecidable sets below $A^\jump$. 

\section The supertask machines

Let us quickly review here how the machines work. For a full account, see \cite[HamLew]. The basic idea is to extend the operation of an ordinary Turing machine into transfinite ordinal time by defining the limit stages of computation. For convenience, we use a three-tape Turing machine model, with separate input, scratch, and output tapes. The successor stages of computation proceed, just as for an ordinary Turing machine, according to a finite program running with finitely many states. Thus, the head reads the tape, reflects on its state and then writes $0$ or $1$ on the tape and moves left, right or not at all according to the program's rigid instructions. The new behavior occurs at a limit stage: the head is reset to the initial starting position; the machine is placed in the special {\it limit} state, just another of the finitely many states; and the values on the cells of the tape are updated by computing the $\limsup$ of the previous cell values. With the limit stage configuration thus completely specified, the machine simply continues the computation. If the {\it halt} state is eventually reached, the machine gives as output whatever is written on the output tape. Since there is plenty of time for the machines to deal with infinite input and output, the natural context for the machines is Cantor space $2^\w$, which in this paper we denote by $\R$ and refer to as the set of reals. Thus, the machines provide notions of computable partial functions $f\from\R\to\R$ on the reals, as well as notions of decidable and semi-decidable sets of reals $A\of\R$. 

Just as with Turing machines, we can use a real as an oracle by simply adding an extra oracle tape on which the oracle real is written, and allow the machine freely to consult this oracle tape during the computations. But because the machines provide a notion of decidability and semi-decidability for {\it sets} of reals $A\of\R$, we really want an oracle concept also in this case, though of course we can't expect to somehow write out an uncountable set $A$ on the tape. Instead, for an oracle $A\of\R$, we allow the machine to make membership queries of $A$, by writing out any real $x$ on a special oracle tape and when in the {\it query} state receiving the answer {\it yes} or {\it no} whether $x\in A$ or not. Thus, in parallel with the notion of constructibility from a predicate, as in the definition of $L[A]$, the machine is allowed to know whether $x\in A$ or not for any $x$ which it is able to produce. The oracle concept leads in turn to the notions of relative decidability and the infinite time Turing degrees. In \cite[HamLew] we presented the basic structural properties of the infinite time Turing degrees, and we aim now to build on that foundation. 

We will follow as closely as possible the notation of \cite[HamLew]. In particular, we write $\varphi_p^A(x)$ to denote the output of program $p$ on input $x$ with oracle $A$, and we write $\varphi_{p,\a}(x)\halts$ to mean that the computation of program $p$ on input $x$ converges in fewer than $\a$ many steps. We will generally denote reals by $a$, $b$, $c$ and $x$, $y$, $z$ and sets of reals by $A$, $B$ and $C$. Thus, when we mention oracles $A$ and $a$ without further explanation, we intend that $A$ should be taken as a set of reals and $a$ as an individual real. Since we proved in \cite[HamLew] that every real $x$ is equivalent to the set $A_x$ of truncated initial segments of $x$ (with $0$s appended), the case of set oracles $A$ subsumes that of real oracles $a$. We write $A\ilte B$ to mean that the characteristic function of $A$ is computable from $B$. This leads naturally to the strict computability relation $\ilt$ and the degree relation $\iequiv$. The lightface halting problem is $$0^\jump=h=\hbox{$\set{p\st\varphi_p(0)\halts}$},$$ with the corresponding lightface jump defined by $A^\jump=A\oplus h^A$ where $h^A$ is the set $\set{p\st \varphi_p^A(0)\halts}$. Similarly, the boldface halting problem is $$0^\Jump=\set{\<p,x>\st \varphi_p(x)\halts},$$ with the corresponding boldface jump defined by $A^\Jump=\set{\<p,x>\st \varphi_p^A(x)\halts}$. In \cite[HamLew] we proved, for example, that $A\ilt A^\jump\ilt A^\Jump$. If $y$ is a real coding a relation $\tlt$ on (a subset of) $\w$ with order type $\a$, then $y\restrict\b$ is the real coding the restriction of the relation $\tlt$ to its first $\b$ many elements. 
An ordinal is {\it writable} if there is a machine which on input $0$ gives as output a real coding that ordinal. An ordinal 
is {\it clockable} when it is the length of a halting infinite time Turing machine computation. As always, we denote the supremum of the writable ordinals by $\l$ and  the supremum of the clockable ordinals by $\g$. A real is {\df eventually} writable when it eventually appears on the output tape during the computation of an infinite time Turing machine on input $0$, never subsequently to be changed (but the machine need not halt). A real is {\df accidently}
writable when it appears on one of the tapes during a computation on input $0$. We say that an ordinal is eventually 
or accidentally writable when it is coded by a real which is, respectively, eventually or accidentally writable. We denote the supremum of the eventually writable ordinals by $\z$ and the supremum of the accidently writable ordinals by $\Sigma$. The results of \cite[HamLew] establish, for example, that $\l\leq\g$ and $\l<\z<\Sigma$. 

The proof of our Main Theorem has been considerably simplified by the recent result of Philip Welch \cite[Wel98a] which answers what was probably the primary open question of \cite[HamLew], namely, the question whether every clockable ordinal is writable. Let us therefore paraphrase Welch's argument here. 

\theorem.(Welch) Every clockable ordinal is writable. In particular, $\l=\g$.
\ref\Welch

\proof By results of \cite[HamLew], it suffices to show that $\z>\g$. (The basic reason for this is that there can be no clockable ordinals between $\l$ and $\z$, indeed, between $\l$ and $\Sigma$, for if $\b$ is clockable and accidently writable, then we can search for a real coding an ordinal which is long enough to support that computation, and halt when we find one. So $\b$ would be writable.  Thus, in order to prove $\g=\l$ it suffices to show that $\z>\g$.)

We begin by showing that in any computation on input $0$, if a cell is $0$ at $\Sigma$, then in fact it must be $0$ for the entire time between $\z$ and $\Sigma$. Thus in fact it must stabilize before $\z$, being $0$ from some point on. Suppose some cell is $0$ at $\Sigma$ for the computation of program $p$ on input $0$. We will describe an algorithm which eventually writes a real coding an ordinal which is at least as large as the stabilizing stage of that cell, the stage at which it turns $0$ never to change again before $\Sigma$. On the side, our algorithm simulates all programs on input $0$, and thereby produces a steady stream of all accidentally writable reals coding ordinals. For each such real $w$ produced, the algorithm simultates the computation of $p$ along the ordinal $\a$ coded by $w$, and looks for the stage $\b$ when the values of the cell have stabilized at $0$, if indeed they have stabilized before $\a$. Using $w$ the algorithm is able to write a real coding $\b$. Now, the algorithm checks the output tape to see if the real on the output tape codes an ordinal at least as large as $\b$. If so, the algorithm leaves the output tape unchanged. If $\b$ is larger than what appears on the output tape, then the algorithm replaces the output tape with the real coding $\b$. The key idea is that eventually the algorithm will see accidentally writable ordinals which are larger than the true stabilizing point below $\Sigma$, and at such a stage, we will write an ordinal which is at least as large as the true stabilizing point on the output tape. After this, the output tape will never change, because no other accidentally writable real can support a computation with a larger stabilizing point than the true stabilizing point. Thus, our algorithm eventually writes an ordinal which is at least as large as the stabilizing point of the computation, so this stabilizing point must be eventually writable, as we claimed. 

Similarly, if a cell is eventually $1$ before $\Sigma$, then in fact it must be $1$ the entire time between $\z$ and $\Sigma$, and in fact must be $1$ constantly from some point before $\z$. What the discussion until now shows is that the stage $\Sigma$ snapshot of any computation on input $0$ is the same as the stage $\z$ snapshot, and this is true in the strong sense that any cell which is $0$ at $\Sigma$ is $0$ between $\z$ and $\Sigma$. Consequently, any computation which does not halt by $\z$ will be endlessly repeating. Consequently, $\z$ is at least as large as $\g$.\qed

The proof proceeded by establishing the following fact.

\lemma.(Welch) Any non-halting supertask computation on input $0$ attains it's repeating snapshot by stage $\z$, and the stage $\z$ snapshot repeats at stage $\Sigma$. 
\ref\Repeat

Of course, we mean this in the strong sense, so that no cell which is $0$ at $\z$ ever turns to $1$ between $\z$ and $\Sigma$; thus, the computation is truly in an endlessly repeating loop. The theorem generalizes easily to any oracle, so that we have $\l^A=\g^A$ for any oracle $A$.

\section A negative solution to Post's problem for supertasks

Let us now prove the main theorem in the context of degrees represented by a real, namely, that there are no such degrees between $0$ and $0^\jump$. Thus, in the reals there are only two semi-decidable supertask degrees. 

\theorem Main Theorem. There are no reals $z$ such that $0\ilt z\ilt 0^\jump$.
\ref\Main

\proof Suppose towards a contradiction that $0\ilt z\ilt 0^\jump$. Thus, with oracle $0^\jump$ there is a program $p$ yielding a halting program on input $0$ with output $z$. Consider now the algorithm with no oracle which computes approximations to $0^\jump$ by simulating all computations and keeping track of which programs have halted. Thus, at any stage, this algorithm has a current best guess for $0^\jump$, namely, the set of programs which have halted on input $0$ by that stage. Eventually, of course, by stage $\g$, the program will have simulated the computations for long enough so that all the programs which will ever halt have already halted, and the {\it true approximation} will be reached, the stage when the best approximation to $0^\jump$ is actually $0^\jump$ itself. Of course, our algorithm has no way of realizing that it has reached the true approximation, and will forever continue to simulate the computations, waiting for more programs to halt. In any case, for each of the approximations to $0^\jump$ along the way, let our algorithm use it as an oracle with program $p$ in an attempt to write the real $z$, and to continue doing so until the next better approximation to $0^\jump$ is found. Eventually, of course, the true approximation to $0^\jump$ will be reached, and this true approximation will support a computation by the program $p$ yielding the real $z$ (while our algorithm continues searching for better approximations to $0^\jump$). This algorithm shows that the real $z$ is eventually writable. Let us think a bit more about this algorithm. Suppose first that $z$ is not produced by any of the proper approximations to $0^\jump$, and first appears only after the true approximation $0^\jump$ has been produced. In this case, with $z$ as an oracle, we could recognize $0^\jump$ as the first approximation to produce $z$ with program $p$. This contradicts our assumption that $z\ilt 0^\jump$. Therefore, second, it must be that $z$ is produced by one of the proper approximations to $0^\jump$, appearing at some stage before $\g$. Since the approximations remain unchanged until the next program halts, the proper approximations are all writable, and consequently $z$ is computable from a writable oracle. Thus, $z$ must itself be writable, contradicting our assumption that $0\ilt z$. We conclude that no such $z$ exists.\qed

The theorem generalizes easily to oracles:

\theorem. For any real $a$ there is no real $z$ such that $a\ilt z\ilt a^\jump$. Indeed, for any oracle $A$ there is no real $z$ such that $A\ilt A\oplus z\ilt A^\jump$. 

\proof Suppose towards a contradiction that $A\ilt A\oplus z\ilt A^\jump=A\oplus h^A$. Thus, with $A$ and $h^A$ as oracles there is a program $p$ yielding a halting computation on input $0$ with output $z$. With only the oracle $A$, now, consider the algorithm which computes approximations to $h^A$ and for each of these approximations, together with $A$, attempts to write $z$ with program $p$. If none of the proper approximations is able to produce $z$ in this way then $z$ will only appear after the true approximation to $h^A$ appears, and so from $A$ and $z$ we could recognize $h^A$, contradicting $A\oplus z\ilt A^\jump$. Thus, it must be that with oracle $A$ and one of the proper approximations to $h^A$, we can write $z$. Since all such proper approximations are $A$-writable, it follows that $z$ is $A$-writable, contradicting $A\ilt A\oplus z$.\qed

The non-existence of real degrees between a real $a$ and it's jump $a^\jump$ means of course that the structure of the supertask degrees differs remarkably from the classical analogue. The Sack's Density theorem, for example, obviously cannot hold in the context of the supertask degrees on the reals. Does it hold for supertask degrees in general, that is, in the context of sets of reals? Let us emphasize that while the previous theorems show that there are no noncomputable {\it reals} below $0^\jump$, the results of section four will show that there is a rich degree structure between $0$ and $0^\jump$ when it comes to {\it sets} of reals. In particular, we will show that there are many degrees $A$ such that $0\ilt A\ilt 0^\jump$. This is what makes the next corollary so intriguing. 

\corollary Lowness Corollary. Every degree below $0^\jump$ is low. That is, if $A\ilt 0^\jump$, then $A^\jump\iequiv 0^\jump$. In fact, for any degree $A$, if $0^\jump\not\ilte A$ then $A^\jump=A\oplus 0^\jump$. 

\proof Certainly $A\oplus0^\jump\ilte A^\jump$ since $A^\jump$ computes both $A$ and $0^\jump$. If $0^\jump\not\ilte A$ and $A\oplus 0^\jump\ilt A^\jump$, then $A\ilt A\oplus 0^\jump\ilt A^\jump$, contradicting the previous theorem. Consequently, if $0^\jump\not\ilte A$ then it must be that $A^\jump\iequiv A\oplus 0^\jump$. Therefore, in particular, if $A\ilt 0^\jump$ then certainly $0^\jump\not\ilte A$ and so $A^\jump\iequiv A\oplus 0^\jump\iequiv 0^\jump$. So $A$ is low, as desired. \qed

\corollary. More generally, for any real $b\ilte A$, if $b^\jump\not\ilte A$, then $A^\jump\iequiv A\oplus b^\jump$. And if $b\ilte A\ilt b^\jump$ then $A^\jump\iequiv b^\jump$. 

\proof This general argument is no more difficult. Suppose that $b\ilte A$ and $b^\jump\not\ilte A$. Then certainly $A\oplus b^\jump\ilte A^\jump$. If they are not equivalent, then $A\ilt A\oplus b^\jump\ilt A^\jump$, contradicting the theorem. Similarly, if $b\ilte A\ilt b^\jump$ then $A^\jump\iequiv A\oplus b^\jump\iequiv b^\jump$, as desired.\qed

Perhaps the main idea behind the proof of the Main Theorem is the following principle: 

\theorem. For any oracle $A$ we have $0^\jump\ilte A$ if and only if $\l<\l^A$. Indeed, whenever $A\ilte B$, then $A^\jump\ilte B$ if and only if $\l^A<\l^B$. 

\proof If $B$ computes $A^\jump$ then $B$ can write a real coding the pre-order in which $p\tlt q$ when $\varphi_p^A(0)$ halts before $\varphi_q^A(0)$, which by results in \cite[HamLew] has order-type $\l^A$. Conversely, if $B$ can write a real coding $\l^A$, then $B$ knows how to long to wait for the $A$-computations to halt, and so it can compute $h^A$. Since $A\ilte B$, this means that $A^\jump\ilte B$.\qed

\section The eventually and accidentally writable degrees

Given that in the reals there are no degrees between $0$ and $0^\jump$, and more generally no reals between any real $a$ and $a^\jump$, it is natural to iterate the jump $\jump$ operator and the build a kind of backbone to the supertask degrees, namely, the following sequence $$0,0^\jump,0^{\jump^2},\ldots,0^{\jump^\a},\ldots$$
We simply apply the $\jump$ operator at successor stages, and take a kind of sum at limit stages, using a real coding the ordinal in question. Specifically, we inductively define $0^{\jump^\a_y}$, whenever $y$ is a real coding a relation $\tlt$ on $\w$ with order type $\a$. At successor stages, if $y$ codes $\a+1$, we simply apply the $\jump$ operator: 
$$0^{\jump^{\a+1}_y}=(0^{\jump^\a_{y\restrict\a}})^\jump.$$ 
At limit stages, when $y$ codes a limit ordinal $\a$, we take a sum:
$$0^{\jump^\a_y}=\oplus_y\set{0^{\jump^\b_{y\restrict\b}}\st\b<\a},$$
where we use the real $y$ to organize the information. That is, with suitable coding $0^{\jump^\a_y}$ is the set $\set{\<n,0^{\jump^\b_{y\restrict\b}}>\st\hbox{$n$ has order type $\b$ in $\tlt$}}$. Of course, in the notation $0^{\jump^\a_y}$ the superscript $\a$ is redundant, since $y$ codes $\a$, but we nevertheless retain it for the suggestion that we are iterating the $\jump$ operator; later for eventually writable ordinals $\a$ we will give an independent meaning to $0^{\jump^\a}$ free of any code for $\a$. Perhaps the clearest way to picture $0^{\jump^\a_y}$ is as an $\a\cross\w$ matrix whose rows are organized using the relation coded by $y$; the $0^\th$ row is trivial, the next row is always the $\jump$ of the previous row, and the limit rows are the sum, organized using the corresponding restriction of $y$, of the previous rows. It is a relatively straightforward exercise, by considering the programs which halt if a specific bit of the oracle is on or off, to verify that $y$ itself is computable from $0^{\jump^\a_y}$. 

There is, of course, a natural connection between this backbone sequence and the eventually writable degrees. In \cite[HamLew] and in the proof of Theorem \Main, we argued that the halting problem $0^\jump$ is eventually writable; one simply embarks on the supertask algorithm which simulates all programs on input $0$, waits for them to halt and keeps track of which ones have halted. At every stage one has an approximation to $0^\jump$, namely, the set of programs which have halted by that stage. Eventually, of course, the simulations will have proceeded long enough that every program which will ever halt has already halted, and the {\it true stage} is reached, the stage by which the approximation to $0^\jump$ is $0^\jump$ itself. The conclusion is that $0^\jump$ is eventually writable. 

The same argument shows that if a real $z$ is eventually writable, then so is $z^\jump$; one simply embarks on the computation which eventually writes $z$, and for each approximation to $z$ encountered along the way, embark on the simulation of all computations which on input $0$ use this real as an oracle and, by keeping track of which programs have halted, obtain an approximation to $z^\jump$. Eventually, the true approximation to $z$ appears, and after this, after all the computations which halt using $z$ as an oracle have halted, the correct $z^\jump$ is eventually written. For any real $z$ let us define the {\df proper approximations} to $z^\jump$ as the approximations which appear on the output tape, using $z$ as an oracle, before the true approximation appears. Thus, the proper approximations to $z^\jump$ are simply the sets of programs which, with $z$ as an oracle, halt before some given $z$-clockable stage.

More generally still, piling approximation upon approximation, one can iterate this idea given a code $y$ for any ordinal $\a$. One simply uses $y$ to organize an $\a\cross\w$ matrix, placing all $0$s in the $0^\th$ row, and using each row as an oracle to approximate its jump in the next row. Each limit row $\b$ is simply the sum of the information in the previous rows, using $y\restrict\b$ to organize the sum. Thus, for example, the algorithm waits for programs to halt in order to approximate $0^\jump$ in the first row, and for each such approximation it waits for programs to halt on that oracle to approximate $0^{\jump^2}$ in the second row, and so on. Each time a new program halts in some row $\b$, we update the approximation in row $\b$ and completely erase the approximations appearing in all rows above $\b$ (since they depended on incorrect information in row $\b$). In order to avoid the accumulating garbage at limits of these erasures, we also completely erase any row at a stage which is a limit of stages where we erased that row (since these approximations also depend on incorrect information). Eventually, the first row will attain the true approximation to $0^\jump$, and will never subsequently be erased. Gradually every row in turn will similarly stabilize, and in the end we will have eventually written $0^{\jump^\a_y}$. What the procedure shows is that $0^{\jump^\a_y}$ is eventually $y$-writable. Consequently, if $y$ itself is eventually writable, then $0^{\jump^\a_y}$ is eventually writable. Perhaps this procedure, in which the approximations to an eventually writable real are used to generate further approximations which are in turn used for further approximations and so on transfinitely, can best be viewed as a kind of transfinite-injury priority construction. We injure the approximations to the later iterates when we learn new, better information about the earlier iterates. 

Let us define the proper approximations to $0^{\jump^\a_y}$, with oracle $y$ coding $\a$, as the reals coding the $\a\cross\w$ matrices (organized using $y$) which appear in the preceding transfinite-injury computation before the true approximation appears. 

\lemma Proper Approximation Lemma. The proper approximations to $z^\jump$ made with oracle $z$ are all $z$-writable. Similarly, if $y$ codes $\a$, then the proper approximations to $0^{\jump^\a_y}$ are $y\oplus 0^{\jump^\b_{y\restrict\b}}$-writable for some $\b<\a$. 

\proof The proper approximations to $z^\jump$ are simply the sets of programs which, using oracle $z$, halt by a given $z$-clockable stage. Such sets are clearly $z$-writable. 

Suppose that $z$ is the proper approximation to $0^{\jump^\a_y}$ appearing at stage $\xi$ in the transfinite-injury computation just previously described. By stage $\xi$ there is a certain maximal portion $0^{\jump^\b_{y\restrict\b}}$ of the oracle which has been correctly computed, and the approximations are only changing in the rows from $\b$ to $\a$. Since the approximation is not yet correct in row $\b$, there must be some stage $\d>\xi$ and a program $p$ which halts in row $\b$, in the approximation to $0^{\jump^{\b+1}_{y\restrict\b+1}}$. And since this is the first time for $p$ to halt in row $\b$ after the true $0^{\jump^\b_{y\restrict\b}}$ has appeared, it must be that $\d$ is $y\oplus 0^{\jump^\b_{y\restrict\b}}$-clockable. Consequently, $\xi$ is $y\oplus 0^{\jump^\b_{y\restrict\b}}$-writable. Thus, with $y\oplus 0^{\jump^\b_{y\restrict\b}}$ as an oracle, we can stop the original construction at stage $\xi$, and thereby write $z$, as desired.\qed

The next lemma is a generalization of the Main Theorem. Whereas the Main Theorem shows that there are no reals between any $a$ and $a^\jump$, the following lemma shows more, that no reals sneak into the limit stages in the transfinite iteration of the $\jump$ operator. 

\lemma Continuity Lemma. If $y$ is an $x$-writable real coding the ordinal $\a$ and $x\ilte 0^{\jump^\a_y}$, then $x\iequiv y\oplus 0^{\jump^\b_{y\restrict\b}}$ for some $\b\leq\a$. 

\proof Let us suppress the subscript $y$, using always $y$ or the appropriate restriction $y\restrict\b$ as the case may be. Consider the algorithm which with oracle $x$ first writes $y$ and then computes approximations to $0^{\jump^\a}$ by recursion along the order given by $y$, in the transfinite-injury manner previously described. Since $x$ is computed by some program $p$ on input $0$ using the oracle $0^{\jump^\a}$, we can for each of our approximations simulate the operation of program $p$ using the approximation as our oracle. If $x$ only appears after the true approximation $0^{\jump^\a}$ is produced, then we could recognize $0^{\jump^\a}$ as the first approximation able to produce $x$ and thereby conclude $x\iequiv y\oplus 0^{\jump^\a}$. Alternatively, suppose one of the proper approximations to $0^{\jump^\a}$ is able to produce $x$. By the previous lemma, all such proper approximations are $y\oplus 0^{\jump^\b}$-writable for some $\b<\a$. And since we could use $x$ to recognize the first approximation to produce $x$, we conclude that $x\iequiv y\oplus 0^{\jump^\b}$ for some $\b<\a$, as desired.\qed

\corollary Continuity Corollary. In particular, if $y$ is a writable real coding $\a$, then every $x$ computable from $0^{\jump^\a_y}$ is equivalent to $0^{\jump^\b_{y\restrict\b}}$ for some $\b\leq\a$.
\ref\ContinuityCor

Often, we only care about the degree of $0^{\jump^\a_y}$, rather than the particular instantiation with the particular code $y$; the next lemma shows that provided one does not try to sneak undue complexity into these codes, the degree of $0^{\jump^\a_y}$ does not depend on $y$. Specifically, we define that $y$ and $z$, both coding the same ordinal $\a$, are {\df hereditarily equivalent} when for every $\b\leq\a$ we have $y\restrict\b\iequiv z\restrict\b$. Any two writable reals coding the same ordinal, for example, are hereditarily equivalent. 

\lemma. If $y$ and $z$ are hereditarily equivalent, then $0^{\jump^\a_y}\iequiv 0^{\jump^\a_z}$. 

\proof We prove this by a straightforward induction on $\a$, the ordinal coded by $y$ and $z$. For the successor case, suppose that $y$ and $z$ both code the ordinal $\a+1$ and inductively, that $0^{\jump^\a_{y\restrict\a}}\iequiv 0^{\jump^\a_{z\restrict\a}}$. By taking the $\jump$ of both sides, it follows that $0^{\jump^{\a+1}_y}\iequiv 0^{\jump^{\a+1}_z}$. For the limit case, suppose that $y$ and $z$ both code the limit ordinal $\a$, and that $y\restrict\b\iequiv z\restrict\b$ for all $\b\leq\a$. Inductively, we may assume that $0^{\jump^\b_{y\restrict\b}}\iequiv 0^{\jump^\b_{z\restrict\b}}$ for $\b<\a$. Since at $\a$ we simply take the sums of these, we need only argue by symmetry that from $0^{\jump^\a_y}$ we can uniformly compute $0^{\jump^\b_{z\restrict\b}}$. Given $0^{\jump^\a_y}$ we can write $y$ and $z$ and begin to compute each $0^{\jump^\b_{z\restrict\b}}$. At limit rows we simply take the sum of what we have computed so far. At successor rows, given $0^{\jump^\b_{z\restrict\b}}$, we search for (and find) the program which witnesses the forward reduction in the equivalence $0^{\jump^\b_{y\restrict\b}}\iequiv 0^{\jump^\b_{z\restrict\b}}$ and then use that program to translate the jump of the former (which we can compute from $0^{\jump^\a_y}$) to the jump of the latter, thereby producing $0^{\jump^{\b+1}_{z\restrict\b+1}}$, as desired.\qed

We now define $0^{\jump^\a}$, when $\a$ is a writable ordinal, to be the degree of $0^{\jump^\a_y}$ where $y$ is a writable real coding $\a$. Since any two writable reals coding $\a$ are hereditarily equivalent, $0^{\jump^\a}$ is well-defined. We thereby obtain the backbone sequence of degrees $$0,0^\jump,\ldots,0^{\jump^\a},\ldots$$ by iterating through the writable ordinals.

\theorem Backbone Continuity Theorem. The backbone sequence, iterated through the writable ordinals, is continuous; namely, whenever $\a$ is a writable limit ordinal, $0^{\jump^\a}$ is a minimal upper bound of $0^{\jump^\b}$ for $\b<\a$.

\proof This follows immediately by Corollary \ContinuityCor.\qed

By the end of this section, we aim to show the same result for the iterations through the eventually writable ordinals and what is more, that this backbone sequence exhausts the eventually writable degrees. In particular, the eventually writable degrees are well-ordered! Thus, in the previous theorem we will actually know that $0^{\jump^\a}$ is the least
eventually writable upper bound to the $0^{\jump^\b}$. In order to prove these things, let us develop a bit of the theory of the eventually and accidentally writable degrees.  We say that a real {\it appears} at time $\a$ when it is written on one of the tapes at stage $\a$ during the (possibly non-halting) supertask computation of some program $p$ on input $0$.

\theorem Timing Theorem.(Welch) Just as the writable reals are those appearing before time $\l$ on an infinite time Turing computation with input $0$, the eventually writable reals are those which appear before time $\z$, and the accidentally writable reals are those which appear before time $\Sigma$.
\ref\Timing

\proof Perhaps this is implicit in \cite[Wel98a]. Every writable real appears of course at a clockable stage, and therefore before stage $\g=\l$. Conversely, it is not difficult to argue that any real which appears at a writable stage is writable. Skipping now to the accidentally writable reals, Lemma \Repeat\ shows that any computation on input $0$ either halts before $\l$ or else obtains its repeating snapshot by stage $\z$, repeating it by stage $\Sigma$. Thus, any accidentally writable real must appear before stage $\Sigma$. And obviously any real which appears before $\Sigma$ is accidentally writable. Considering now the eventually writable reals, if a real $z$ appears before stage $\z$ during the computation of program $p$ on input $0$, then we could eventually write a real coding the ordinal stage when it appeared, and simultaneously eventually write the real appearing at that stage in the computation of $p$ on input $0$, thereby eventually writing $z$. Conversely, if a real $z$ is eventually writable, then since it must appear before stage $\Sigma$ we could, while computing approximations to $z$, search for a real coding an ordinal long enough to see that $z$ has appeared, and thereby succeed in eventually writing a real coding the ordinal stage at which $z$ first appeared. So this must be an eventually writable ordinal.\qed

\lemma Longest Computation Lemma. There is a single program $p$ which on input $0$ takes as long as is possible to repeat, namely, it first attains its repeating snapshot at stage $\z$, first again repeated at stage $\Sigma$. 

\proof This lemma was also observed by Welch. Note that Lemma \Repeat\ establishes that the stage $\z$ snapshot of any supertask computation on input $0$ repeats again at stage $\Sigma$. The issue is whether this snapshot appears also before $\z$ or between $\z$ and $\Sigma$. 

We claim first that the eventually writable reals cannot all appear by some stage $\b$ bounded below $\z$. To see this for any particular $\b$, consider the algorithm which eventually writes a real coding $\b$, and for each approximation $\b'$ to $\b$ makes a table with $\b'$ many rows and simulates all programs on input $0$ for $\b'$ many steps of computation. After this, the algorithm diagonalizes against this table to write a real which does not appear on any computation before stage $\b'$, and then continues to search for better approximations to $\b$. When the true approximation to $\b$ is achieved, this algorithm will eventually write a real which does not appear before stage $\b$, as we claimed. 

Let now $p$ be the program which on input $0$ simultaneously simulates all computations on input $0$. By the observation of the previous paragraph, this algorithm will not reach its repeating snapshot before stage $\z$. Since Lemma \Repeat\ shows that the stage $\z$ snapshot repeats at stage $\Sigma$, it remains only to show that the stage $\z$ snapshot does not repeat before stage $\Sigma$. Consider now the algorithm which searches for all accidentally writable reals which code ordinals, and tests them to see if they code an ordinal long enough to reach the repeating stage of our program $p$. If $p$ repeats at some stage $\a<\Sigma$, then such a real will eventually be found, and this second algorithm could write a real coding $\a$. But certainly $\a$ can not be writable, since $p$ does not repeat until $\z$, so the program $p$ must not repeat again until stage $\Sigma$.\qed

What the argument establishes is that any computation which does not attain its repeating snapshot until stage $\z$, the latest possible stage, will not repeat this snapshot again until stage $\Sigma$, the longest possible delay.

\theorem Linearity Theorem. The accidentally writable reals are linearly ordered by $\ilte$, and this order agrees with their order of first appearance. Namely, if $x$ appears accidentally before or at the same time as $y$, then $x\ilte y$. Hence, the accidentally writable reals are actually well-ordered by $\ilte$. 
\ref\Linear

\proof Of course, we really mean pre-orders here, because one must work modulo $\iequiv$ in order to have an actual order. Suppose that $x$ and $y$ are accidentally writable, appearing first at stages $\a$ and $\b$ by programs $p$ and $q$, respectively, on input $0$, and that $\a\leq\b$. By the Timing Theorem \Timing, it must be that $\b$ is below $\Sigma$. Thus, using $y$ as an oracle, we can search for a real $z$ coding an ordinal which is large enough for $y$ to appear in the simulation of $q$, that is, which is at least $\b$. Since we could easily output such a real $z$ when we found it, the first such $z$ we encounter is computable from $y$. (And actually, from $z$ we could run the computation on which $y$ appears for exactly $\b$ steps, thereby concluding $y\iequiv z$.) Furthermore, there is some natural number $n$ which is the $\a^\th$ element in the order coded by $z$, and so using $z$ we can run the computation on which $x$ appears for exactly $\a$ many steps, and thereby write $x$. So $x\ilte z\ilte y$, as desired.\qed

Let us say that a degree is accidentally or eventually writable, respectively, when there is a representative real in the degree which is accidentally or eventually writable. 

\corollary. The eventually writable degrees have order type $\z$. 

\proof If a real $y$ is computable from an eventually writable real $z$, then $y$ is also eventually writable, and furthermore, the entire computation sequence showing that $y\ilte z$ is eventually writable. Thus, for any eventually writable ordinal $\b$, we can search for and eventually write a listing of (the first representatives encountered for) the first $\b$ many eventually writable reals. Thus, since this list is eventually writable, it does not contain all the eventually writable reals. Thus, there are at least $\z$ many eventually writable degrees. There can be no more than this on account of Theorems \Linear\ and \Timing.\qed

Now it is clear how we may continue the backbone sequence of degrees into the eventually writable ordinals. At successor stages we simply apply the $\jump$, so that $$0^{\jump^{\a+1}}=(0^{\jump^\a})^\jump.$$ And at limit stages $\a$, since the eventually writable degrees are well-ordered, we simply take a supremum $$0^{\jump^\a}=\sup_{\b<\a}0^{\jump^\b}.$$ Since this definition agrees with the earlier one when $\a$ is writable, we have a uniform definition of the backbone sequence of degrees $$0,0^\jump,\ldots,0^{\jump^\a}\!\!,\ldots$$ where we iterate through the eventually writable ordinals. And by construction, using the order topology on the ordinals and on the eventually writable degrees, this sequence is continuous.

\theorem Backbone Theorem. The eventually writable degrees are exactly the degrees on the backbone sequence.

\proof Certainly all the degrees on the backbone sequence are eventually writable. Conversely, by the Main Theorem no degrees appear between a backbone degree $0^{\jump^\a}$ and its successor $0^{\jump^{\a+1}}\!$; by continuity no degrees can sneak in at the limit stages; and by the previous corollary no eventually writable degrees appear on top of or to the side of the backbone sequence. Thus, the backbone sequence exhausts all the eventually writable degrees.\qed

Although we have defined $0^{\jump^\a}$ as the supremum of the previous degrees when $\a$ is an eventually writable limit ordinal, we would rather define every $0^{\jump^\a}$ concretely as the degree of $0^{\jump^\a_y}$ for some eventually writable real $y$ coding $\a$, in analogy with our definition for the writable ordinals. The problem with this approach is that even if one always uses a $y$ of least possible degree, two such $y$ need not necessarily be hereditarily equivalent; indeed, under this definition $0^{\jump^\a}$ would be ill-defined, because some such $y$ could code a lot of extra information into $y\restrict\w$, for example, and then $0^{\jump^\w_{y\restrict\w}}$ would be much more complex than $0^{\jump^\w}\!$. And the same could happen at larger eventually writable ordinals. One way to salvage this problem is to use smooth codes: an eventually writable real $y$ coding an ordinal $\a$ is {\df smooth} when for all $\b\leq\a$ the real $y\restrict\b$ has least possible degree among reals coding $\b$. Thus, the smooth codes in a sense have hereditarily least degree coding the ordinal in question, and sneak as little information into the limit stages as possible. One can easily show that $0^{\jump^\a_y}\iequiv 0^{\jump^\a}$ whenever $y$ is a smooth code for $\a$, and so this could provide an alternative definition of $0^{\jump^\a}\!$, concretely connected with $\jump$-iterations. But does every eventually writable ordinal have a smooth code? 

We conclude this section by showing how to cap off the backbone sequence with the accidentally writable degrees. The fact is that we need add only one more degree! This was also observed by Welch \cite[Wel98b]. 

\theorem Maximal Degree Theorem. The accidentally writable reals have order type $\z+1$ under $\ilt$. In particular, there is a unique accidentally writable degree which is not eventually writable. This degree is therefore maximal among all accidentally writable degrees.

\proof If suffices to argue that there is only one accidentally writable degree which is not eventually writable. Let $p$ be the program as above which on input $0$ simulates all programs on input $0$, and let $z$ be the stage $\z$ snapshot of this computation. We argued before that the snapshot $z$ appears first at stage $\z$ and then is first repeated again at stage $\Sigma$. Therefore, and this is the key idea, $\Sigma$ is $z$-clockable, since with $z$ as an oracle we could simply wait for the second appearance of the snapshot $z$. It follows by the relativized version of Theorem \Welch\ that $\Sigma$ is $z$-writable. Suppose now that $x$ is another accidentally writable real, appearing first at some stage $\a$. If $\a<\z$, then $x$ is eventually writable (and so $x\ilte z$ since $z$ appears first at stage $\z$). Otherwise, $\z\leq\a<\Sigma$. Since $z$ appears at $\z$, we know that $z\ilte x$, and so it suffices to argue that $x\ilte z$. Since $\Sigma$ is $z$-writable, it follows also that $\a$ and $\a+1$ are also $z$-writable, and so with $z$ as an oracle, we can write a real coding $\a+1$, and then using this real we can simulate the computation on which $x$ appears at stage $\a$, and thereby write $x$. So $x\ilte z$.\qed

It is natural therefore to denote this maximal accidentally writable degree by $0^{\jump^\z}\!$, since in the accidentally writable degrees it is the supremum of the $0^{\jump^\a}$ for $\a<\z$. The following list $$0,0^\jump,\ldots,0^{\jump^\a}\!\!,\ldots,0^{\jump^\z}$$ is therefore a complete list of all the accidentally writable degrees. One naturally wants to extend the sequence further, computing $0^{\jump^{\z+1}}$ and so. How far can one go? 

\remark. The accidentally writable reals are not closed under $\iequiv$. 

\proof The proof of the previous theorem shows that $\Sigma$ is $z$-writable for the maximal accidentally writable real $z$. That is, there is a real $y\ilte z$ which codes $\Sigma$. Since $z$ appears at stage $\Sigma$ of a computation, it follows also that $z\ilte y$ and so $z\iequiv y$. But since no accidentally writable real can code $\Sigma$, it must be that $y$ is not accidentally writable.\qed

Some of the ideas and results of this section are considerably improved upon by Philip Welch in \cite[Wel98b], where he shows for example that the degrees $0^{\jump^\a}$ are precisely the master codes of certain models $L_{\l^\a}$. We therefore refer the reader to that paper for many interesting further theorems. The following question remains open: 

\question. In the context of the reals, is $0^{\jump^\w}$ the least upper bound of the sequence $0$, $0^\jump$, $0^{\jump\jump}$, $\ldots$? 

We have shown above that $0^{\jump^\w}$ is the least upper bound of the $0^{\jump^n}$ among the set of accidentally writable reals. And since any real computable from $0^{\jump^\w}$ is eventually writable, it follows more generally that $0^{\jump^\w}$ is a minimal upper bound of the $0^{\jump^n}$ among the set of all real degrees (and more generally that $0^{\jump^\b}$ is a minimal upper bound of $0^{\jump^\a}$ for $\a<\b$ whenever $\b$ is an eventually writable ordinal). What remains is the question of whether there are other upper bounds to this sequence which are incomparable to $0^{\jump^\w}\!$. Techniques in \cite[Wel98c] establish that, say, a real which is Cohen generic over the $V$ will be incomparable with $0^\jump$, and since the existence of such incomparable degrees is a $\Sigma^1_2$ fact, it must by absoluteness be true already in the ground model $V$. Thus, there are numerous reals incomparable to the $0^{\jump^n}\!$. Nevertheless, the above question remains open.

\section A positive solution to Post's problem for supertasks

We now prove the second part of the Main Theorem, that in the context of sets of reals, there are degrees strictly between $0$ and $0^\jump$. 

\theorem Main Theorem. There are incomparable semi-decidable sets of reals. Since such sets are properly between $0$ and $0^\jump$, the supertask analogue of Post's Problem has a positive solution. 
\ref\MainTheorem

We will adapt the classical priority argument of Friedburg-Munchnik to the supertask context. This is not a completely straightforward enterprise, for two reasons. First, in the classical argument one can preserve an oracle computation by preserving finitely much information about that oracle and promising not to injure this information later; but in the supertask context an oracle computation can use infinitely much information about the oracle, and we may have difficulty preserving this much information. We will get around this difficulty with a strategic application of the Boundedness Lemma, below, by only adding writable reals to the oracle which have not appeared as queries on the previous computations. The second, more serious difficulty is that in the classical context oracle computations can be diagonalized against by simply handling finite length computations using approximations to the oracles; but in the supertask context an oracle can potentially yield very long halting computations which cannot be foreseen by an algorithm equipped only with approximations to that oracle. Simply put, oracle computations can easily halt beyond $\g$, and we may have difficulty diagonalizing against such computations. Our solution to this difficulty will lie in the Reflection Lemma, below. We will take care to ensure that the oracles that we construct do not yield new clockable ordinals above $\g$. Indeed, this precaution is absolutely necessary, because if an oracle can halt beyond $\g$, then it would be able to compute $0^\jump$ and therefore could not possibly be the kind of oracle we desire.

\lemma Boundedness Lemma. There is no halting computation on input $0$ in which every writable real accidently appears. Hence, for any clockable ordinal $\a$ there is a writable real which does not appear accidentally on any computation before $\a$.

\proof Since $\a$ is clockable, it is also writable, and consequently the entire sequence of snapshots detailing the computation of length $\a$ is writable. Thus, all these snapshots are also writable, and by diagonalizing against the reals appearing on the computation we can output a writable real not on the computation. By carrying out this procedure with respect to all programs simultaneously, we obtain a single writable real which does not appear on any computation on input $0$ before stage $\a$.\qed

Our general technique for constructing a semi-decidable oracle $A$ will be to gradually place writable reals into $A$ as the construction proceeds. The oracle $A$, therefore, will be the union of the approximation oracles $A_\a$, the set of reals added before stage $\a$. We promise to only add reals to $A$ at clockable stages, and because of this the {\it true approximation} $A_\g$, which appears at stage $\g$, will be $A$ itself. We imagine that $A$ is represented on the tape as the reals which appear as rows in an $\w\cross\w$ matrix (suitably coded on a portion of the scratch tape). If at some clockable stage we want to add to $A$ the writable real $x$, produced by program $p$, then we will write $x$ in the $p^\th$ row of this matrix. In particular, a writable real will only be added to $A$ after the program which produced this real has halted. We imagine that in the background of the construction we are simulating all programs on input $0$, thereby producing a steady stream of writable reals and information about which programs have produced which reals and in which order they halted. The approximation $A_\a$ is simply the set of reals which have appeared as rows in the matrix by stage $\a$. Since no supertask algorithm on input $0$ can halt beyond $\g$, it is impossible for our procedure to know when it has reached the true stage, and so even when the true approximation $A=A_\g$ has been written out on the tape, at stage $\g$, the algorithm will nevertheless continue trying to build better approximations to $A$, not realizing that it has exhausted the clockable ordinals. This interesting state of affairs, in which the true oracle has been correctly written on the tape, but the algorithm does not know that it has finished, provides the key technique of our argument. 

\lemma Reflection Lemma. Suppose a set $A$, consisting entirely of writable reals, is the semi-decidable union of approximations $A_\a$ according to a computable procedure in the manner described above, and that $\varphi_p^A(x)\halts=0$ for some writable real $x$. Then there are unboundedly many clockable ordinals $\a<\g$ such that $\varphi_{p,\a}^{A_\a}(x)\halts=0$. 

\proof Consider the algorithm which searches for a nonclockable stage $\b$ such that the approximation $A_\b$ supports a halting computation $\varphi_p^{A_\b}(x)\halts=0$, and furthermore that this computation halts before the next clockable ordinal (i.e. before the next better approximation to $A$). This algorithm must succeed in its search, since if it doesn't succeed before $\g$, then the true approximation $A_\g$ supports a halting computation $\varphi_p^A(x)=0$. But since no algorithm can halt beyond $\g$, this algorithm in fact must succeed before $\g$. Thus, there is a least clockable ordinal $\a$ beyond the $\b$ found by the algorithm. Since reals are only added at clockable stages, $A_\b=A_\a$, and the computation $\varphi_p^{A_\a}(x)=\varphi_p^{A_\b}(x)=0$ halts before $\a$, as desired. An easy modification shows in addition that such $\a$ must occur unboundedly often.\qed

Now we are ready to prove the second Main Theorem.

\proofof Main Theorem \MainTheorem: We will adapt the classical Freidburg-Munchnik argument to the supertask context. For the reasons explained earlier, we must take some care. We will construct two sets $A$ and $B$, consisting entirely of writable reals, by specifying in a supertask computable way their approximations $A_\a$ and $B_\a$ as described above. We want to meet the requirements $$R_p:\varphi_p^B\not=A,\qquad\qquad S_p:\varphi_p^A\not=B.$$ The earlier requirements have higher priority, and we will only injure a requirement to satisfy a higher priority requirement. Each requirement will be injured only finitely often. At every stage we will associate with each requirement $R_p$ and $S_p$ a respective witness $x$, and this witness will be a real which is writable but which does not appear even accidentally on any computation of length $r(p,\a)$, a clockable ordinal which we call the restraint (represented for the purposes of this algorithm by a program which clocks it). For any clockable ordinal, such reals exist by the Boundedness Lemma. If it happens at a clockable stage $\a$ that $\varphi_{p,\a}^{B_\a}(x)\halts=0$, where $x$ is the witness currently attached to $R_p$, then we say that the requirement $R_p$ {\it requires attention}. And now it {\it receives attention} when we add $x$ to $A_{\a+1}$ and redefine the restraint function $r(p,\a+1)$ to preserve, with priority $p$, the computation $\varphi_p^B(x)$. That is, we don't want in the future to add any elements to $B$ which appeared as queries on this computation. Since $B_\a$ can be coded with a writable real---our algorithm has written $B_\a$ on a portion of the tape by stage $\a$---it follows that the computation $\varphi^{B_\a}_p(x)$ could be simulated by a computation without any oracle, and therefore by the Boundedness Lemma there are writable reals which do not appear even accidently during this computation. In particular, they do not appear as queries to $B$. Thus, for each of the weaker priority requirements, we search for new witnesses which do not appear as a query to $B$ on the computation $\varphi_p^{B_\a}(x)$ (or on such computatons for the stronger priority requirements). And we are guarranteed to find them. With priority $p$ this will ensure that $\varphi_p^B(x)=\varphi_p^{B_\a}(x)=0$. Since $x\in A$, we thereby satisfy requirment $R_p$. We perform an analogous procedure to give attention to the requirements $S_p$ when they require it, by adding the witnesses to $B$. 

It is clear that if a requirement requires attention and then receives it, and this attention is not injured later by any stronger priority requirment, then the requirement is satisfied. Suppose, alternatively, that the requirement $R_p$ does not require attention after some stage. This means that for sufficiently large clockable $\a$, the computation $\varphi_{p,\a}^{B_\a}(x)$ does not converge to $0$. Also, the witness $x$ was never added to $A$ (because there are only finitely many strong priority requirements, from some point on we did not change the witness $x$). If $\varphi_p^B(x)\halts=0$ then by the Reflection Lemma, it must be that for unboundedly many clockable $\a$, $\varphi_{p,\a}^{B_\a}(x)\halts=0$, a contradiction. So either $\varphi_p^B(x)$ does not converge at all or it gives some other output than $0$. Therefore, since $x\notin A$, the requirement $R_p$ is satisfied. A similar argument shows $S_p$ is satisfied, and so the construction satisfies all the requirments. Consequently, neither $A$ nor $B$ is computable from the other, and the theorem is proved.\qed Main Theorem

Naturally, many open questions remain. What is the nature of the supertask computability degrees? To what extent do they resemble the Turing degrees?

\section Bibliography

\nopagenumbers
\parindent=0pt
\newbox\Article
\newbox\Journal
\newbox\Author
\newbox\Vol
\newbox\No
\newbox\Year
\newbox\Page
\newbox\Book
\newbox\Publisher
\newbox\Pubaddr
\newbox\Key
\newbox\Editor
\newbox\Comment
\newbox\Note
\def\entry#1#2\par{\item{#1\quad}\hskip-1.1em#2\par}
\def\article#1{\setbox\Article=\hbox{\sl #1, }}
\def\journal#1{\setbox\Journal=\hbox{\rm #1 }}
\def\author#1{\setbox\Author=\hbox{\sc #1, }}
\def\vol#1{\setbox\Vol=\hbox{\bf #1 }}
\def\no#1{\setbox\No=\hbox{no. #1 }}
\def\year#1{\setbox\Year=\hbox{\rm({\oldstyle #1}) }}
\def\page#1{\setbox\Page=\hbox{\rm p. #1 }}
\def\book#1{\setbox\Book=\hbox{\it #1, }}
\def\publisher#1{\setbox\Publisher=\hbox{\rm #1, }}
\def\pubaddr#1{\setbox\Pubaddr=\hbox{\rm #1, }}
\def\key#1{\setbox\Key=\hbox{#1}}
\def\editor#1{\setbox\Editor=\hbox{\rm(#1, Ed.) }}
\def\comment#1{\setbox\Comment=\hbox{\rm #1}}
\def\note#1{\setbox\Note=\hbox{\rm #1 }}
\def\ref#1\par{\smallskip{#1
  \entry{\ifhbox\Key\unhbox\Key\else[\ ]\fi}%
  \unhbox\Author\unhbox\Note
  \ifhbox\Book \unhbox\Book\unhbox\Publisher\unhbox\Pubaddr
               \unhbox\Editor\unhbox\Page\unhbox\Year\unhbox\Comment
  \else \unhbox\Article\unhbox\Journal\unhbox\Vol\unhbox\No\unhbox\Editor
        \unhbox\Page\unhbox\Year\unhbox\Comment\fi\par}}

\tenpoint

\ref
\author{Joel David Hamkins and Andrew Lewis}
\article{Infinite Time Turing Machines}
\journal{to appear in the Journal of Symbolic Logic}
\key{[HamLew]}

\ref
\author{Robert Soare}
\book{Recursively Enumerable Sets and Degrees}
\publisher{Springer-Verlag}
\year{1980}
\key{[Soa]}

\ref
\author{Philip Welch}
\article{The lengths of infinite time Turing machine computations}
\journal{submitted to Bulletin and Journal of the London Math. Soc.}
\key{[Wel98a]}

\ref
\author{Philip Welch}
\article{Eventually infinite time Turing degrees}
\comment{(in preparation)}
\key{[Wel98b]}

\ref
\author{Philip Welch}
\article{Friedman's trick: minimality in the infinite time Turing degrees}
\journal{submitted to ``Sets and Proofs', vol. 2, 1997 Proc. European Meeting of the ASL, Leeds, (CVP LMS Lecture Notes in Maths Series 1998)}
\key{[Wel98c]}

\bye